\documentclass[graybox]{svmult}

\usepackage{type1cm}        
\usepackage{makeidx}         
\usepackage{graphicx}        
\usepackage{multicol}        
\usepackage[bottom]{footmisc}

\usepackage{newtxtext}       %
\usepackage[varvw]{newtxmath}       

\makeindex             

\def\a{\mathbf{a}}
\def\b{\mathbf{b}}
\def\d{\mathbf{d}}
\def\e{\mathbf{e}}
\def\f{\mathbf{f}}
\def\g{\mathbf{g}}
\def\h{\mathbf{h}}
\def\n{\mathbf{n}}
\def\q{\mathbf{q}}
\def\w{\mathbf{w}}
\def\x{\mathbf{x}}
\def\w{\mathbf{w}}
\def\z{\mathbf{z}}
\def\zero{\mathbf{0}}

\def\dt{\partial_t}
\def\ddt{\frac{d}{dt}}
\def\dx{\partial_x}

\def\E{\mathcal{E}}
\def\H{\mathcal{H}}
\def\T{\mathcal{T}}

\def\curl{\operatorname{curl}}

\usepackage{tikz}
\usepackage{pgfplots}
\pgfplotsset{compat=newest}

\title*{On higher order passivity preserving schemes for nonlinear Maxwell's equations}

\author{Herbert Egger and  Vsevolod Shashkov}
\institute{Herbert Egger \at Johann Radon Institute for Computational and Applied Mathematics and Johannes Kepler Universität Linz, Altenbergerstr. 69, 4040 Linz, Austria \email{herbert.egger@jku.at}
\and Vsevolod Shashkov \at TU Darmstadt, Dolivostr. 15, 64293 Darmstadt, Germany \email{shashkov@mathematik.tu-darmstadt.de}}

\begin{document}
\maketitle
\abstract{
We present two strategies for designing passivity preserving higher order discretization methods for Maxwell's equations in nonlinear Kerr-type media. 
Both approaches are based on variational approximation schemes in space and time. This allows to rigorously prove energy conservation or dissipation, and thus passivity, on the fully discrete level. 
For linear media, the proposed methods coincide with certain combinations of mixed finite element and implicit Runge-Kutta schemes. The order optimal convergence rates, which can thus be expected for linear problems, are also observed for nonlinear problems in the numerical tests. 
}

\section{Introduction}
\label{sec:1}

We consider the propagation of high-intensity electromagnetic waves through dielectric media which is of relevance in typical applications of nonlinear optics \cite{Boyd2008,Shen1994}.
The underlying physics are modelled by Maxwell's equations
\begin{equation}
    \dt\d = \curl \h, \qquad \dt \b = -\curl \e \label{eq:general_maxwell}
\end{equation}
with $\e$ and $\h$ denoting the electric and magnetic field intensities and $\d$ and $\b$ the corresponding fluxes.
For the following discussion, we assume that the relation between fields and fluxes is given by
\begin{equation}
    \d = \d(\e):=\epsilon_0(\chi^{(1)} + \chi^{(3)} |\e|^2)\e, \qquad \b = \b(\h):= \mu_0\h \label{eq:material_properties}
\end{equation}
which describes the instantaneous electric and magnetic response of a nonlinear Kerr-type medium.
Here $\epsilon_0$, $\mu_0$ represent the permittivity and permeability of vacuum and the constants $\chi^{(1)}$ and $\chi^{(3)}$ describe the linear and nonlinear dielectric effects. 
Let us note that more general nonlinear constitutive equations as well as lossy materials can be considered with slight modifications of our arguments. 

\medskip 

\textbf{Energy balance.}
The constitutive equations \eqref{eq:material_properties} allow to eliminate the fluxes $\d,\b$ in equation \eqref{eq:general_maxwell} and to rewrite the system solely in terms of field intensities $\e,\h$ as
\begin{align}
    \d'(\e)\dt\e &= \curl\h, \label{eq:eh_form_1}\\
    \mu_0\dt\h &= -\curl \e, \label{eq:eh_form_2} 
\end{align}
where $\d'(\e)=\epsilon_0(\chi^{(1)} + 3\chi^{(3)} |\e|^2)$ denotes the incremental permittivity. 
For the further discussion, we restrict our considerations to a bounded domain $\Omega$ and we assume homogeneous magnetic boundary conditions
\begin{align} \label{eq:eh_form_3}
\n \times \h &= \zero \qquad \text{on } \partial\Omega.
\end{align}
For prescribed fields $\e$, $\h$, the electromagnetic energy of the system is then given by
\begin{equation}
    \E (\e, \h) = \int_\Omega
    w_E(\e) + w_M(\h) \ d\x
\end{equation}
with energy densities 
$w_E(\e) = \frac{\epsilon_0}{2} \left(\chi^{(1)}|\e|^2 +\frac{3\chi^{(3)}}{2}|\e|^4 \right)$
and $w_M(\h)=\frac{\mu_0}{2} |\h|^2$. 
Let us note that these are convex functions of their arguments and further observe that $w_E'(\e) = \d'(\e) \e$ and $w'_M(\h)=\mu_0 \h = \b'(\h) \h$. 
This allows to show that 
\begin{align} \label{eq:energy_balance}
    \E (\e(t), \h(t)) = \E (\e(s),\h(s)) 
\end{align}
for any solution of $(\e,\h)$ of \eqref{eq:eh_form_1}--\eqref{eq:eh_form_3} and all $t \ge s$; see Section~\ref{sec:2}.
Hence the energy of the system is conserved for all time and, in particular, the system is passive. 

\medskip

\textbf{Passivity preserving discretization.}
The finite-difference time-domain method is certainly the industry standard for simulation of Maxwell's equations in time domain. Various extensions to nonlinear media have been proposed, e.g., in \cite{FujiiEtAl2004,JosephTaflove1997,SpachmannEtAl1999,ziolkowski1993full}.  We refer to \cite{MaksymovEtAl2011} for comparison of different approaches and to  \cite{BokilEtAl2018,JiaLiFangLi2019} for recent developments. 
More flexible finite volume and related discontinuous Galerkin approximations have been discussed in \cite{aregba2014godunov} and \cite{huang2017second}. 
A combination of mixed finite element or discontinuous Galerkin approximations with Runge-Kutta time-stepping schemes has been studied in \cite{AanesAngermann2020,FisherEtAl2007} and \cite{blank2013discontinuous,bokil2017energy}. 
While these methods are capable of providing higher order approximations, their strict passivity for higher order approximations in time seems not completely settled. 

\medskip

\textbf{Contributions.}
In this paper, we propose two strategies that allow to construct provably passivity preserving discretization schemes of arbitrary high order. A key ingredient here is the use of a variational time integration methodology \cite{akrivis2011galerkin,LewEtAl2004}. 
As illustrated in  \cite{egger2019structure,egger2021energy}, this allows to obtain conservative or dissipative numerical approximations for a large class of evolution problems.
Our first approach is based on a mixed finite element approximation of the $\e$--$\h$ formulation \eqref{eq:eh_form_1}--\eqref{eq:eh_form_3} in space and a discontinuous Galerkin method in time. The resulting scheme is slightly dissipative and leads to a discrete analogue of \eqref{eq:energy_balance} with an inequality instead of equality. 
The second scheme is based on a reformulation of the problem in terms of the electric field $\e$ and the magnetic vector potential $\a$.
We consider a discretization of this $\e$--$\a$ formulation in space with the same finite element spaces for both variables and a Petrov-Galerkin time integration. This leads to a scheme with exact energy conservation. 

\medskip

\textbf{Outline.}
In Section~\ref{sec:2}, we derive the variational form of \eqref{eq:eh_form_1}--\eqref{eq:eh_form_3} and prove the energy balance \eqref{eq:energy_balance}. 
Section~\ref{sec:3} then discusses the discretization of the $\e$--$\h$ formulation and presents the corresponding discrete energy inequality.
In Section~\ref{sec:4}, we derive the equivalent $\e$--$\a$ formulation, its variational formulation, and the corresponding energy balance. 
Section~\ref{sec:5} is devoted to the discretization of this problem for which we prove exact energy  conservation on the discrete level. 
In Section~\ref{sec:6}, we illustrate our theoretical findings by numerical tests and comment on the convergence properties and performance of the proposed methods. 
We then close with a short summary.

\section{The $\e$--$\h$ formulation}
\label{sec:2}

Let us briefly introduce the most relevant notation used in the rest of the paper. 
We write $L^2(\Omega)$ for the space of square integrable scalar or vector valued functions and denote by $H(\curl;\Omega) = \{\w \in L^2(\Omega) : \curl \w \in L^2(\Omega)\}$ the set of functions with square integrable derivatives. 
Furthermore, we use  $\langle \f,\g\rangle = \int_\Omega \f \cdot \g \, d\x$ to abbreviate the $L^2$-scalar product. The following lemma summarizes the main properties of solution to the $\e$-$\h$ formulation that we will use in the following.
\begin{lemma} 
Let $(\e,\h)$ be a smooth solution of \eqref{eq:eh_form_1}--\eqref{eq:eh_form_3}. Then 
\begin{align}
    \langle\d'(\e(t))\dt\e(t), \w\rangle &= \langle \h(t), \curl \w\rangle \label{eq:eh_var_1}\\
    \langle\mu_0 \dt \h(t), \q \rangle &= -\langle \curl \e(t), \q\rangle \label{eq:eh_var_2}
\end{align}
for all $\w \in H(\curl;\Omega)$, $\q \in L^2(\Omega)$, and all $t \ge 0$. Moreover, 
\begin{align} \label{eq:eh_energy_identity}
    \E(\e(t),\h(t)) &= \E(\e(s),\h(s)) \qquad \forall t \ge s.
\end{align}
\end{lemma}
\begin{proof}
The second equation follows immediately by multiplying \eqref{eq:eh_form_2} with the test function $\q$ and integration over the domain $\Omega$. To verify the first identity, we  multiply equation \eqref{eq:eh_form_1} by $\w$, integrate over $\Omega$, and observe that 
\begin{align*}
\langle \d'(\e) \dt \e, \w \rangle
&= \langle \curl \h,\w\rangle
 = \langle \h,\curl \w\rangle + \int_{\partial\Omega} \n \times \h \cdot \w \, ds(\x).
\end{align*}
Here we used integration-by-parts in the second step. 
Due to \eqref{eq:eh_form_3}, the boundary term vanishes and we obtain \eqref{eq:eh_var_1}.
By formal differentiation of the energy, we further get
\begin{align*}
\ddt \E(\e,\h) 
&=\langle w_E'(\e), \dt \e\rangle + \langle w_M'(\h), \dt \h\rangle
 = \langle \d'(\e) \e, \dt \e\rangle + \langle \mu_0 \h, \dt h\rangle \\
&= \langle \d'(\e) \dt \e, \e\rangle + \langle \mu_0 \dt \h, \h\rangle 
 = \langle \h, \curl \e\rangle - \langle \curl \e, \h\rangle = 0. 
\end{align*}
In the first line, we used the definition of the energy functional and the relations between the energy densities and the material parameters stated in the introduction.
In the second line, we first exchanged the order of the functions in the scalar product and then used the variational identities \eqref{eq:eh_var_1}--\eqref{eq:eh_var_2} with test functions $\w=\e$ and $\q=\h$. 
The energy identity then follows immediately by integration over time.
\end{proof}
Let us note that in the essential steps of the energy identity, we only made use of the variational form \eqref{eq:eh_var_1}--\eqref{eq:eh_var_2} of the problem with test functions $\w=\e(t)$ and $\q=\h(t)$. This motivates to consider a variational discretization scheme in the following.

\section{Discretization of the $\e$--$\h$ formulation}
\label{sec:3}

Let $W_h \subset H(\curl;\Omega)$ and $Q_h \subset L^2(\Omega)$ denote some finite dimensional subspaces and let $I_\tau=\{t^n: 0\le n \le N \}$ be a sequence of discrete time steps $t^n=n \tau$ with $\tau=T/N$. 
We write $I^n=[t^{n-1},t^n]$ for the $n$th time interval and denote by 
$P_k(I^n;X)$ the space of polynomial functions $v : I_n \to X$ with values in some vector space $X$.  By $(*)|_{t^n}$ we mean the evaluation of the time dependent expression $(*)$ at $t=t^{n}$. 
For discretization of problem \eqref{eq:eh_form_1}--\eqref{eq:eh_form_3}, we then consider the following method.

\begin{problem} \label{prob:eh_discrete}
Let the discrete initial values $\e^0_h \in W_h$ and $\h^0_h \in Q_h$ be given. \\
Then for $1 \le n \le N$ find $\e^n_h \in P_k(I^n;W_h)$ and $\h_h^n \in P_k(I^n;Q_h)$  such that 
\begin{align}
\int_{I_n} \langle \d'(\e^n_h) \dt \e^n_h,  \w_h\rangle -\langle \h^n, \curl  \w_h\rangle dt 
&= 
\langle \d'(\e_h^n)(\e_h^{n-1}-\e_h^{n}), \w_h\rangle|_{t^{n-1}}  \label{eq:eh_var_1h}\\
\int_{I_n} \langle \mu_0 \dt \h^n_h,  \q_h\rangle dt + \langle \curl e_h^n, \q_h\rangle 
&= \langle \mu_0 (\h_h^{n-1}-\h_h^n), \q_h\rangle|_{t^{n-1}} \label{eq:eh_var_2h}
\end{align}
holds for all test functions $\w_h \in P_k(I^n;W_h)$ and $\q_h \in P_k(I^n;Q_h)$. 
\end{problem}

This scheme is based on a Galerkin approximation of \eqref{eq:eh_var_1}--\eqref{eq:eh_var_2} in space and a discontinuous Galerkin method in time \cite{akrivis2011galerkin}. It emerges as a particular example of an abstract discretization framework for dissipative evolution problems; see  \cite[Sec.~8]{egger2019structure}. 
From the theoretical results derived in this reference, we conclude the following. 
\begin{lemma} \label{lem:eh_discrete}
Let $(\e_h^n,\h_h^n)_n$ denote a solution of Problem~\ref{prob:eh_discrete}. Then 
\begin{align} \label{eq:eh_discrete}
    \E(\e_h^n(t^n),\h_h^n(t^n)) \le \E(\e_h^{m}(t^m),\h_h^m(t^m)) \qquad \forall m \le n.
\end{align}
\end{lemma}
\begin{proof}
For $m=n-1$, the inequality above can be verified by testing the variational principle \eqref{eq:eh_var_1h}--\eqref{eq:eh_var_2h} with $\w_h=\e_h^n$ and $\q_h = \h_h^n$ and using convexity of the energy densities $w_E(\e)$ and $w_M(\h)$. The case $m < n-1$ can then be treated by induction. For details of the arguments used in the proof, we refer to \cite[Sec.~3 and 4]{egger2019structure}.
\end{proof}

Let us emphasize that the variational form of the time integration scheme is important here to prove passivity in the case of nonlinear constitutive equations.

\section{The $\e$--$\a$ formulation}
\label{sec:4}

We now present an alternative approach towards the passivity preserving discretization of the problem, which is based on a standard reformulation in terms of the magnetic vector potential. 
Let us start with an auxiliary observation.
\begin{lemma}
Let $(\e,\h)$ be smooth functions satisfying identity \eqref{eq:eh_form_2}. 
Further define $\a(t) = \a_0-\int_0^t \e(s) ds$ 
with $\a_0$ chosen such that $\curl \a_0 = \mu_0 \h(0)$. Then  
\begin{align*}
    \e(t) = - \dt \a(t) 
    \qquad \text{and} \qquad 
    \mu_0 \h(t) = \curl \a(t) 
    \qquad \forall t \ge 0. 
\end{align*}
\end{lemma} 
\begin{proof}
The first identity is clear. Using that $\mu_0$ is constant and employing \eqref{eq:eh_form_2}, we deduce that 
$\mu_0 \dt \h = -\curl \e  = \curl \dt \a.
$
By integration in time, we then get
\begin{align*}
\mu_0 \h(t) 
&= \mu_0 \h(0) + \int_0^t \mu_0 \dt \h \, dt
 = \curl \a_0 + \int_0^t \curl \dt \a  \, dt
 = \curl \a(t),
\end{align*}
which already proves the second identity of the lemma.
\end{proof}
Let us note that $\curl \a = \mu_0 \h = \b$, i.e., $\a$ is just the usual magnetic vector potential frequently used in the magneto-quasi static setting. 
Using the above observations, we can now reformulate the system \eqref{eq:eh_form_1}--\eqref{eq:eh_form_2} equivalently as 
\begin{align}
-\d'(\e) \dt \a &= \d'(\e) \e \label{eq:ea_form_1}\\
\d'(\e) \dt \e &= \curl (\nu_0 \curl \a) \label{eq:ea_form_2}
\end{align}
where we introduced $\nu_0=\mu_0^{-1}$ for convenience. The particular choice of the multiplying factors in the first equation will become clear from the proof of Lemma~\ref{lem:ea} below. 
The boundary condition \eqref{eq:eh_form_3} further translates to
\begin{align} \label{eq:ea_form_3}
    \n \times \curl \a &= 0 \qquad \text{on } \partial\Omega.
\end{align}
For obvious reasons, we call \eqref{eq:ea_form_1}--\eqref{eq:ea_form_3} the $\e$--$\a$ formulation of our problem. 
As a final step, we also rewrite the energy functional in terms of the fields $(\e,\a)$, i.e., 
\begin{align}
    \H(\e,\a) 
    = \E(\e,\nu_0 \curl(a)) 
    = \int_\Omega w_E(\e) + \tfrac{\nu_0}{2} |\curl \a|^2  d\x.
\end{align}
Let us again summarize the basic properties of this alternative formulation. 
\begin{lemma} \label{lem:ea}
Let $(\e,\a)$ denote a sufficiently smooth solution of \eqref{eq:ea_form_1}--\eqref{eq:ea_form_3}. Then 
\begin{align}
    -\langle\d'(\e(t))\dt\a(t), \w\rangle &= \langle\d'(\e(t))\e(t),\w\rangle \label{eq:ea_var_1}\\
    \langle\d'(\e(t))\dt \e(t) , \z\rangle &= \langle\nu_0 \curl\a(t), \curl\z\rangle \label{eq:ea_var_2}  
\end{align}
for all test functions $\w,\z\in H(\curl;\Omega)$ and all $t \ge 0$. 
Furthermore 
\begin{align} \label{eq:ea_energy_indentity}
        \H(\e(t),\a(t)) = \H(\e(s),\a(s)) \qquad \forall s \le t.
\end{align}
\end{lemma}
\begin{proof}
The first equation follows immediately by multiplying \eqref{eq:ea_form_1} with $\w$ and integration over the domain. 
In the same manner, we deduce from equation \eqref{eq:ea_form_2} that
\begin{align*}
\langle \d'(\e) \dt \e, \z \rangle 
&=\langle \curl (\nu_0 \curl \a), \z\rangle\\ 
&= \langle \nu_0 \curl \a, \curl \z \rangle + \int_{\partial\Omega} \n \times (\nu_0 \curl \a) \cdot \z \, ds(\x).
\end{align*}
In the second step, we used integration-by-parts. The boundary term vanishes due to the boundary condition \eqref{eq:ea_form_3}, which already leads to \eqref{eq:ea_var_2}.
By formal differentiation of the energy functional and the relation between the electric energy functional $w_E(\e)$ and the constitutive law $\d(\e)$, we can further see that 
\begin{align*}
\ddt\H(\e,\a) 
&= \langle \w_E'(\e), \dt \e\rangle + \langle \nu_0 \curl(\a), \curl \dt \a \rangle \\
&= \langle \d'(\e) \e, \dt \e\rangle + \langle \nu_0 \curl(\a), \curl \dt \a \rangle. 
\end{align*}
This corresponds to the sum of the two terms on the right hand side of \eqref{eq:ea_var_1}--\eqref{eq:ea_var_2} with test functions $\w=\dt \e$ and $\z=\dt \a$. As a consequence, we thus obtain
\begin{align*}
\ddt\H(\e,\a) 
&= -\langle \d'(\e) \dt \a, \dt \e \rangle + \langle \d'(\e) \dt \e, \dt \a\rangle = 0.
\end{align*}
The energy identity \eqref{eq:ea_energy_indentity}
 now follows immediately by integration over time. 
\end{proof}

Note that in the basic step of the proof, we again simply utilized the variational identities \eqref{eq:ea_var_1}--\eqref{eq:ea_var_2} with the particular test functions $\w=\dt \e(t)$ and $\z=\dt \a(t)$. This motivates the variational discretization scheme of the following section.

\section{Discretization of $\e-\a$ formulation} \label{sec:5}

As before, let $W_h\subset H(\curl;\Omega)$ denote some finite dimensional subspace and further recall the notation about the time grid from Section~\ref{sec:3}.
For the numerical approximation of problem \eqref{eq:ea_form_1}--\eqref{eq:ea_form_3}, we then consider the following method. 
\begin{problem} \label{prob:ea_discretization}
Let $\e_h^0,\a_h^0 \in W_h$ be given 
and for $1 \le n \le N$, find $\e_h^n,\a_h^n \in P_{k+1}(I_n;W_h)$ such that $\e_h^n(t^{n-1})=\e_h^{n-1}(t^{n-1})$ and 
$\a_h^n(t^{n-1})=\a_h^{n-1}(t^{n-1})$ as well as 
\begin{alignat}{2}
-\int_{I^n}\langle\d'(\e_h^n)\dt\a_h^n, \tilde \w_h\rangle dt &= \int_{I^n} \langle\d'(\e_h^n)\e_h^n,\tilde \w_h\rangle dt  
\qquad && \forall \tilde \w \in P_k(I^n;W_h)\label{eq:ea_var_1h}\\
\int_{I^n}\langle\d'(\e_h^n)\dt \e_h^n , \tilde \z_h\rangle &= \int_{I^n}\langle\nu_0 \curl\a_h^n, \curl\tilde \z_h\rangle dt \qquad && \forall \tilde \z \in P_k(I^n;W_h). \label{eq:ea_var_2h}  
\end{alignat}
\end{problem}
This scheme is based on a Galerkin approximation of \eqref{eq:ea_var_1}--\eqref{eq:ea_var_2} in space together with a Petrov-Galerkin time discretization \cite{akrivis2011galerkin}. 
Similar methods can be applied for the numerical solution of a wide class of evolution problems; see  \cite{egger2021energy} for examples. 
Let us note that by construction, the discrete solution is continuous in time and can be computed by an implicit time stepping algorithm.
The most important property of the method can be summarized as follows. 
\begin{lemma} \label{lem:ea_discrete}
Let $(\e_h^n,\a_h^n)_{n \ge 0}$ denote a solution of Problem~\ref{prob:ea_discretization}. Then 
\begin{align} \label{eq:ea_discrete}
\H(\e_h^n(t^n),\a_h^n(t^n)) = \H(\e_h^m(t^m),\a_h^m(t^m)) \qquad \forall n \ge m.
\end{align}
\end{lemma}
\begin{proof}
By the fundamental theorem of calculus, we obtain 
\begin{align*}
(*) :&=\H(\e_h^n(t^n),\a_h^n(t^n)) 
- \H(\e_h^n(t^{n-1}),\a_h^n(t^{n-1})) \\
&= \int_{I^n} \ddt \H(\e_h^n,\a_h^n) \, dt 
 = \int_{I^n} \langle w_E'(\e_h^n) , \dt e_h^n\rangle + \langle \nu_0 \curl \a_h^n, \dt \a_h^n\rangle dt.
\end{align*}
Let us recall that $w_E'(\e) = \d'(\e) \e$, which allows us to conclude that 
\begin{align*}
(*) 
&= \int_{I^n} \langle d'(\e_h^n) \e_h^n , \dt e_h^n\rangle + \langle \nu_0 \curl \a_h^n, \dt \a_h^n\rangle dt \\
&= -\int_{I^n} \langle \d'(\e_h^n) \dt \a_h^n, \dt \e_h^n\rangle - \langle \d'(\e_h^n) \dt \e_h^n, \dt \a_h^n\rangle dt = 0.
\end{align*}
In the second step, we here used the variational identities \eqref{eq:ea_var_1h}--\eqref{eq:ea_var_2h} with 
$\tilde \w_h = \dt \e_h^n$ and $\tilde z_h = \dt \a_h^n$, which is admissible by the choice of approximation and test spaces.   
From the conditions for  $\e_h^n(t^{n-1})$ and  $\a_h^n(t^{n-1})$, we can then deduce that 
\begin{align*}
\H(\e_h^n(t^n),\a_h^n(t^n))
&=\H(\e_h^n(t^{n-1}),\a_h^n(t^{n-1}))
 =\H(\e_h^{n-1}(t^{n-1}),\a_h^{n-1}(t^{n-1})).
\end{align*}
This already verifies the energy identity for $m=n-1$. The general case $m<n-1$ can finally again be obtained by induction. 
\end{proof}



\section{Numerical validation}\label{sec:6}

We now illustrate our theoretical results by some numerical tests and comment on the implementation 
of the proposed methods and their convergence behavior.
%

\medskip 

\textbf{Test problem.}
For simplicity, we consider in the sequel a one-dimensional version of problem \eqref{eq:general_maxwell} over the domain $\Omega = (0,1)$. 
In that case $H(\curl;\Omega)=H^1(\Omega)$ and the two $\curl$ operators in \eqref{eq:ea_form_1}--\eqref{eq:ea_form_2} reduce to $\dx$ and $-\dx$, respectively.
The parameters in the material laws \eqref{eq:material_properties} are chosen as $\epsilon_0=\mu_0=\chi^{(1)}=1$ and in order to illustrate the effect of the nonlinear material response, we will consider the two choices $\chi^{(3)}=0$ and $\chi^{(3)}=0.1$ below. 
Note that the problem is linear in the first case. 
The initial values are finally set to $\e(0,x)=\exp(-100 x^2)$ and $\a(0)=\h(0)=\zero$.

\medskip 

\textbf{Spatial approximation. }
Let $\T_h$ be a uniform mesh with grid points $x_i = i h$ and uniform mesh size $h=1/M$.  
We use piecewise polynomial spaces
\begin{align*}
    W_h = P_p(\T_h) \cap H(\curl;\Omega), \qquad Q_h = P_{p-1}(\T_h)
\end{align*}
over the grid $\T_h$ for the space discretization with polynomial degree $p \ge 1$. 

\begin{remark}
In order to facilitate the implementation of the proposed methods, the scalar products $\langle \cdot,\cdot\rangle$ in the discrete variational problems \eqref{eq:eh_var_1h}--\eqref{eq:eh_var_2h} and  \eqref{eq:ea_var_1h}--\eqref{eq:ea_var_2h} are approximated by inexact versions $\langle \cdot,\cdot\rangle_h$, which are realized by numerical quadrature. 
In our computations, we use the Gauss-Lobatto formula with $p+1$ nodes on every element, which is a standard choice; see\cite{CohenMonk1998,cohen2002higher,GeeversEtAl2018}.
Let us note that the discrete energy inequality \eqref{eq:eh_discrete} and the  identity \eqref{eq:ea_discrete} remain valid, if the same quadrature rule is used for defining the discrete energy functionals.  
\end{remark}


\textbf{Time integration.}
Since the material law $\d(\e)$ only involves polynomial nonlinearities, all time integrals in our discretization methods can be computed exactly by  numerical quadrature. 
For the solution of the nonlinear systems in Problems~\ref{prob:eh_discrete} and~\ref{prob:ea_discretization}, we utilize a simple fixed point iteration with tolerance set to $10^{-12}$. 

\medskip 

\begin{figure}[h]
\centering
\begin{tikzpicture}
\footnotesize
\begin{axis}[
width = 0.7\textwidth,
height = 0.25\textwidth,
grid=both, 
minor grid style={gray!25}, 
major grid style={gray!25},
no marks,
width=\linewidth]
\addplot[red, thick] table[x={x}, y={solT0}] {soldata.txt};
\addplot[gray, dashed, thick] table[x={x}, y={solLT0}] {soldata.txt};
\end{axis}
\end{tikzpicture}

\begin{tikzpicture}
\footnotesize
\begin{axis}[
width = 0.7\textwidth,
height = 0.25\textwidth,
grid=both, 
minor grid style={gray!25}, 
major grid style={gray!25},
no marks,
width=\linewidth]
\addplot[red, thick] table[x={x}, y={solT40}] {soldata.txt};
\addplot[gray, dashed, thick] table[x={x}, y={solLT40}] {soldata.txt};
\end{axis}
\end{tikzpicture}

\begin{tikzpicture}
\footnotesize
\begin{axis}[
width = 0.7\textwidth,
height = 0.25\textwidth,
grid=both, 
minor grid style={gray!25}, 
major grid style={gray!25},
no marks,
width=\linewidth]
\addplot[red, thick] table[x={x}, y={solT80}] {soldata.txt};
\addplot[gray, dashed, thick] table[x={x}, y={solLT80}] {soldata.txt};
\end{axis}
\end{tikzpicture}

\begin{tikzpicture}
\footnotesize
\begin{axis}[
width = 0.7\textwidth,
height = 0.25\textwidth,
grid=both, 
minor grid style={gray!25}, 
major grid style={gray!25},
no marks,
width=\linewidth]
\addplot[red, thick] table[x={x}, y={solT120}] {soldata.txt};
\addplot[gray, dashed, thick] table[x={x}, y={solLT120}] {soldata.txt};
\end{axis}
\end{tikzpicture}

\begin{tikzpicture}
\footnotesize
\begin{axis}[
width = 0.7\textwidth,
height = 0.25\textwidth,
grid=both, 
minor grid style={gray!25}, 
major grid style={gray!25},
no marks,
width=\linewidth]
\addplot[red, thick] table[x={x}, y={solT160}] {soldata.txt};
\addplot[gray, dashed, thick] table[x={x}, y={solLT160}] {soldata.txt};
\end{axis}
\end{tikzpicture}

\caption{Snapshot $\e_h^n(t^n)$ of the numerical solution at time steps $t^n=0.0,0.2,0.4,0.6,0.8$ obtained with the method of Section~\ref{sec:5} for two scenarios: linear case ($\chi^{(3)}=0$; black dashed) and nonlinear case ($\chi^{(3)}=0.1$; red solid). In both cases, the discrete energy is preserved exactly.}
\label{fig:wave}
\end{figure}
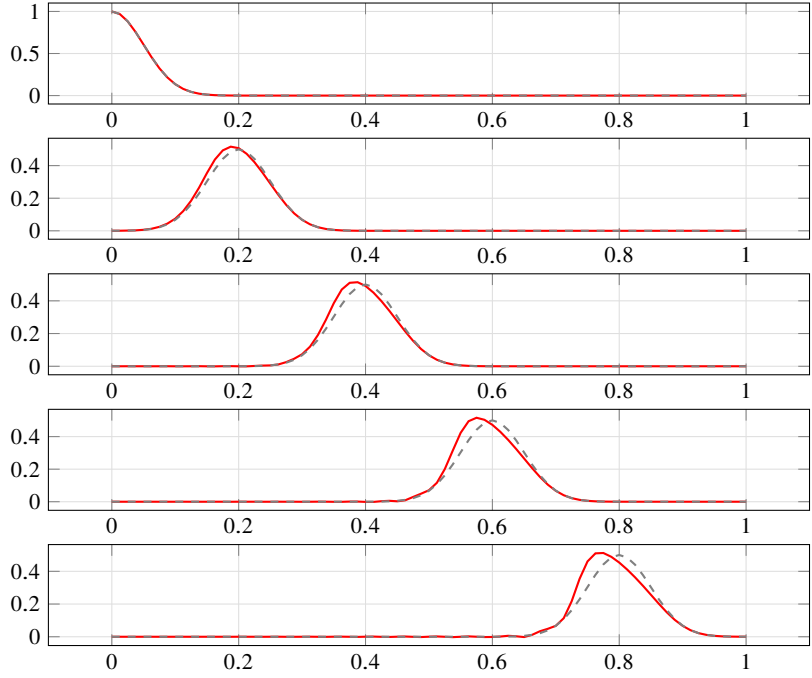

\textbf{Comparison of linear and nonlinear material behavior.}
In Figure~\ref{fig:wave} we display some snapshots $\e_h^n$ of the numerical approximations for the electric field in case of a medium with linear and a nonlinear material behaviour, respectively.

In the absence of nonlinearities,  the signal simply propagates to the right hand side without changing its shape. 
The nonlinear medium response, on the other hand, leads to a slowing down of the propagation at high intensities, which results in a precurser at the beginning and a ramp up at the end of the pulse.  
%
%

\subsection{Results for the $\e$-$\h$ formulation}

Let us note that in the linear case ($\chi^{(3)}=0$), the method \eqref{eq:eh_var_1h}--\eqref{eq:eh_var_2h} corresponds to a mixed finite element approximation in space with time stepping by the Radau-IIA method with $s=k+1$ stages \cite{akrivis2011galerkin}.
From the error analysis in \cite{cohen2002higher,GeeversEtAl2018,Monk1992}, we expect that the error behaves like 
\begin{align}
    err_{h,\tau}^{p,k} := \max_{0 \le t^n \le T} \|\e_h^n(t^n) - \e^{2n}_{h/2}(t^n)\|_{h/2} = O(h^{p+1} + \tau^{2k+1}).
\end{align}
Here $p$ and $k$ are the polynomial degree of the spatial and temporal approximation, respectively, and $\e_h^n$ is the numerical solution on the interval $I_n=[t^{n-1},t^n]$ with $t^n = n \tau$, while we denote by $\e^{2n}_{h/2}$ the corresponding solution obtained on uniformly refined mesh and with time step $\tilde \tau=\tau/2$. 
The approximation $\|\cdot\|_{h/2}$ for the $L^2$-norm is computed by numerical quadrature on the finer grid $\T_{h/2}$. 
In Table~\ref{tab:1} and \ref{tab:2}, we display the errors obtained in our numerical tests for the nonlinear case $\chi^{(3)}=0.1$ for different approximation orders $p$ and $k$ in space and time. 

\begin{table}[ht!]
\centering
\begin{tabular}{c||c|c||c|c||c|c}
$h $ & \multicolumn{2}{c||}{$p=1$} & \multicolumn{2}{c||}{$p=2$} & \multicolumn{2}{c}{$p=3$}\\
\hline
 & \, err $\times \, 10^{-1}$ \, & \, e.o.c. \, & \, err  $\times \, 10^{-2}$ \, & \, e.o.c. \, & \, err $\times \,10^{-3}$ \, & \,e.o.c.\,\\
\hline
0.05    & 0.248186 & ---  & 0.387722 & ---  & 0.417564 & --- \\
0.025   & 0.071272 & 1.80 & 0.003319 & 3.54 & 0.018346 & 4.50 \\
0.0125  & 0.018438 & 1.95 & 0.000299 & 3.47 & 0.000950 & 4.27 \\
0.00625 & 0.004641 & 1.99 & 0.000034 & 3.11 & 0.000058 & 4.02
\end{tabular}
\caption{Convergence of the method of Section~\ref{sec:3} with respect to the mesh size $h$ for different polynomial orders $p$ used for the spatial approximation.\label{tab:1}}
\end{table}

\begin{table}[ht!]
\centering
\begin{tabular}{c||c|c||c|c||c|c}
$\tau $ & \multicolumn{2}{c||}{$k=0$} & \multicolumn{2}{c||}{$k=1$} & \multicolumn{2}{c}{$k=2$}\\
\hline
 & \, err $\times \, 10^{-1}$ \, & \, e.o.c. \, & \, err  $\times \, 10^{-2}$ \, & \, e.o.c. \, & \, err $\times \,10^{-3}$ \, & \,e.o.c.\,\\
\hline
0.025    & 0.257057 & ---  & 0.280420 & ---  & 0.550798 & --- \\
0.0125   & 0.171025 & 0.61 & 0.038358 & 2.87 & 0.019199 & 4.84 \\
0.00625  & 0.100673 & 0.76 & 0.004879 & 2.98 & 0.000610 & 4.97 \\
0.003125 & 0.054697 & 0.88 & 0.000612 & 3.00 & 0.000019 & 5.00
\end{tabular}
\caption{Convergence of the method of Section~\ref{sec:3} with respect to the time step size $\tau$ for different polynomial orders $k$ used for the temporal approximation.\label{tab:2}}
\end{table}



The convergence rates that are expected for the linear case are also observed for the nonlinear case. In all computations, the discrete energy decays monotonically. Since the solutions are uniformly bounded, one can see that, similar to the linear case, the energy error behaves like 
$$
|\E(\e,\h) - \E(\tilde \e,\tilde \h))| \approx \|\e-\tilde \e\|^2+\|\h-\tilde \h\|^2.
$$
As a consequence, we expect and observe very small energy errors. In particular for higher order methods, the numerical dissipation seems therefore negligible.

\subsection{Numerical results for the $\e$--$\a$ formulation}

As a next step, we investigate the convergence of our second discretization scheme.
The spatial approximation here is a standard finite element method and, in the linear case ($\chi^{(3)}=0$), the time discretization amounts to the Lobatto-IIIA method with $s=k+1$ stages. 
We thus expect that the error behaves like 
\begin{align}
    err_{h,\tau}^{p,k} := \max_{0 \le t^n \le T} \|\e_h^n(t^n) - \e^{2n}_{h/2}(t^n)\|_{h/2} = O(h^{p+1} + \tau^{2k+2}).
\end{align}
In Table~\ref{tab:3} and \ref{tab:4}, we  display the numerical errors obtained with the method of Section~\ref{sec:5} for different approximation orders $p$ and $k$ in space and time. 

\begin{table}[ht!]
\centering
\begin{tabular}{c||c|c||c|c||c|c}
$h $ & \multicolumn{2}{c||}{$p=1$} & \multicolumn{2}{c||}{$p=2$} & \multicolumn{2}{c}{$p=3$}\\
\hline
 & \, err $\times \, 10^{-3}$ \, & \, e.o.c. \, & \, err  $\times \, 10^{-3}$ \, & \, e.o.c. \, & \, err $\times \,10^{-3}$ \, & \,e.o.c.\,\\
\hline
0.05    & 0.412735 & ---  & 0.297889 & ---  & 0.277589 & --- \\
0.025   & 0.127333 & 1.70 & 0.022976 & 3.69 & 0.011844 & 4.55 \\
0.0125  & 0.033235 & 1.94 & 0.002874 & 2.99 & 0.000747 & 3.99 \\
0.00625 & 0.008372 & 1.99 & 0.000359 & 3.00 & 0.000046 & 3.99
\end{tabular}
\caption{Convergence of the method of Section~\ref{sec:5} with respect to the mesh size $h$ for different polynomial orders $p$ used for the spatial approximation.\label{tab:3}}
\end{table}

\begin{table}[ht!]
\centering
\begin{tabular}{c||c|c||c|c||c|c}
$\tau $ & \multicolumn{2}{c||}{$k=0$} & \multicolumn{2}{c||}{$k=1$} & \multicolumn{2}{c}{$k=2$}\\
\hline
 & \, err $\times \, 10^{-1}$ \, & \, e.o.c. \, & \, err  $\times \, 10^{-3}$ \, & \, e.o.c. \, & \, err $\times \,10^{-4}$ \, & \,e.o.c.\,\\
\hline
0.05    & 0.801343 & ---  & 0.611080 & ---  & 0.368882 & --- \\
0.025   & 0.226645 & 1.82 & 0.040060 & 3.93 & 0.006549 & 5.81 \\
0.0125  & 0.057709 & 1.97 & 0.002538 & 3.98 & 0.000108 & 5.93 \\
0.00625 & 0.014537 & 1.99 & 0.000160 & 3.98 & 0.000002 & 5.96
\end{tabular}
\caption{Convergence of the method of Section~\ref{sec:5} with respect to the time step size $\tau$ for different polynomial orders $k$ used for the temporal approximation.\label{tab:4}}
\end{table}

As can be easily be deduced from the tables, the convergence rates are again exactly as expected. Let us further mention that the discrete energy was preserved up to round-off errors in all our computational tests with this method.

\section{Discussion}

In this paper, we discussed two different approaches towards the construction of higher order provably passivity preserving numerical schemes for Maxwell's equations in nonlinear media. A key ingredient was the use of appropriate variational space and time discretization schemes which allowed us to rigorously prove fully discrete energy identites, respectively, inequalities on the discrete level. 
Both approaches investigateded in the paper lead to implicit time-stepping schemes, which for linear media coincide with certain Runge-Kutta methods. 
The proposed schemes show the expected convergence behavior for linear as well as for nonlinear problems. A full error analysis should be possible but is left for future research.


\bibliographystyle{spmpsci}
\bibliography{literature}

\end{document}